\documentclass[a4paper,12pt]{article}
\usepackage{latexsym}
\usepackage{graphicx}
\usepackage{amsthm,amsmath,amsfonts,amssymb,graphicx,graphics,psfrag}

\newtheorem{thm}{Theorem}
\newtheorem{oq}{Question}

\newtheorem{prop}{Proposition}

\date{}

\begin{document}

\title{Star-uniform Graphs}

\author{ Mikio Kano$^{1}$, Yunjian Wu$^{2}$ and Qinglin Yu$^{2,3}$
\thanks{Corresponding email: yu@tru.ca}
\\ {\small $^1$ Department of Computer and Information Sciences}
\\ {\small Ibaraki University, Hitachi, Ibaraki, Japan}
\\ {\small $^2$ Center for Combinatorics, LPMC}
\\ {\small Nankai University, Tianjin, China}
\\ {\small $^3$ Department of Mathematics and Statistics}
\\ {\small Thompson Rivers University, Kamloops, BC, Canada}
}

\maketitle

\begin{abstract}
 A {\it star-factor} of a graph $G$ is a spanning subgraph of $G$
such that each of its component is a star.  Clearly, every graph
without isolated vertices has a star factor. A graph $G$ is called
{\it star-uniform} if all star-factors of $G$ have the same number
of components.  To characterize star-uniform graphs was an open
problem posed by Hartnell and Rall, which is motivated by the
minimum cost spanning tree and the optimal assignment problems. We
use the concepts of factor-criticality and domination number to
characterize all star-uniform graphs with the minimum degree at
least two. Our proof is heavily relied on Gallai-Edmonds Matching
Structure Theorem.
\begin{flushleft}
\end{flushleft}
\end{abstract}

\vskip 1mm \vskip 1mm

\section{Introduction}
Throughout this paper, all graphs considered are simple. We refer
the reader to \cite{bb} and \cite{lm} for standard graph theoretic
terms not defined in the paper.

Let $G=(V,E)$ be a graph with vertex set $V(G)$ and edge set $E(G)$.
We call a graph with only one vertex {\it trivial} and all other
graphs {\it nontrivial}. If $S \subset V(G)$, then $G-S$ is the
subgraph of $G$ obtained by deleting the vertices in $S$ and all the
edges incident with them. Similarly, if $E'\subset E(G)$, then
$G-E'=(V(G), E(G)-E')$. The set of vertices adjacent to $S$ in $G$
is denoted by $N_{G}(S)$. If $G$ is not a forest, then the length of
a shortest cycle in $G$ is called its {\it girth}, denoted by
$g(G)$, and a cycle of order $g(G)$ is called a {\it girth cycle}.
An {\it odd $($or even$)$ cycle $($or path$)$} is the one with odd
$($or even$)$ number of vertices.  The {\it union} $G_{1} \cup
G_{2}$ of the graphs $G_{1}$ and $G_{2}$ is the graph with vertex
set $V(G_{1})\cup V(G_{2})$ and edge set $E(G_{1})\cup E(G_{2})$. A
{\it star} is a graph isomorphic to a complete bipartite graph
$K_{1,n}$ for some $n \geq 1$, and the vertex of degree $n$ is
called the {\it center} of the star.

A {\it star-factor} of a graph $G$ is a spanning subgraph of $G$
each component of which is a star. It is not hard to see that every
graph without isolated vertices admits a star-factor. If one limits
the sizes of the stars used, such a star-factor may not always
exist. In \cite{amahashi}, Amahashi and Kano presented a criterion
for the existence of a restricted star-factor, i.e., \{$K_{1, 1},
\cdots , K_{1, n}$\}-factor. Yu \cite{yu} obtained an upper bound on
the maximum number of edges in a graph with a unique star-factor.

  A vertex subset $S$ of a graph $G$ is a {\it dominating set} if
every vertex of $G$ either belongs to $S$ or is adjacent to a vertex
of $S$. The cardinality of a smallest dominating set is called the
{\it domination number} of $G$, denoted by $\gamma(G)$. For
extensive bibliographies regarding work on domination in graphs the
reader is referred to \cite{hay}.

A {\it matching} $M$ of $G$  is a subset of $E(G)$ such that no two
elements of $M$ are adjacent. The number of edges in $M$ is called
the {\em size} of $M$. A matching $M$ is called a {\it maximum
matching} if $G$ has no matching $M'$ with $|M'|>|M|$. If every
$v\in V(G)$ is incident to some edge in $M$, then $M$ is said to be
a {\it perfect matching}. A {\it near-perfect matching} in a graph
$G$ is one covering all but exactly one vertex of $G$. A graph $G$
is said to be {\it factor-critical} if $G-v$ has a perfect matching
for every $v\in V(G)$ and this concept is first introduced by Gallai
\cite{lm}. An $M$-{\it alternating path} $($or {\it $M$-alternating
cycle}$)$ in $G$ is a path $($or cycle$)$ whose edges are
alternately in $M$ and $E-M$. Let $M_{1}$ and $M_{2}$ be matchings
in $G$ and $M_{1} \cup M_{2}$ denote the subgraph formed by the
union of the two edge sets, so $V(M_{1}\cup M_{2}) = V(M_{1})\cup
V(M_{2})$ and $E(M_{1}\cup M_{2}) = E(M_{1})\cup E(M_{2})$. The
components of this subgraph are edges, alternating even cycles or
alternating paths.

Let $G$ be a graph. Denote by $D(G)$ the set of all vertices in $G$
which are not covered by at least one maximum matching of $G$,
$A(G)$ be the set of vertices in $V(G)-D(G)$ adjacent to at least
one vertex in $D(G)$. Finally let $C(G)=V(G)-A(G)-D(G)$ (see Figure
1).

    Gallai \cite{ed} and Edmonds (see \cite{lm}), independently, obtained the
following canonical decomposition theorem for maximum matching in
graphs. This result can be considered as a refinement of Tutte's
famous $1$-Factor Theorem and it provides a complete structural
characterization of the maximum matchings in graphs. It was the
foundation for Edmonds' well-known polynomial algorithm for finding
a maximum matching in graphs. Its full power is still waiting to be
explored in the future.

\begin{thm}[Gallai-Edmonds Structure Theorem for Matchings]\label{gaed}
Let $G$ be a graph and $D(G)$, $A(G)$ and $C(G)$ be the sets defined
as above. Then

$(a)$ every component of the subgraph induced by $D(G)$ is
factor-critical;

$(b)$ the subgraph induced by $C(G)$ has a perfect matching;

$(c)$ the order of $A(G)$ is less than the number of the components
in $D(G)$;

$(d)$ every maximum matching of $G$ consists of a near-perfect
matching of each component of $D(G)$, a perfect matching of each
component of $C(G)$ and a matching which matches all vertices of
$A(G)$ with vertices in distinct components of $D(G)$.
\end{thm}

\begin{figure}[h,t]
\setlength{\unitlength}{1pt}
\begin{center}
\begin{picture}(40,100)(0,60)

\put(-10,150){\line(1,0){90}}\put(-12,148){\line(0,-1){58}}\put(84,148){\line(0,-1){58}}
\put(48,90){\line(1,0){35}}

\put(48,90){\line(-2,-1){57}}\put(-12,90){\line(0,-1){28}}

\put(-12,57){\line(-1,-1){15}}\put(-10,57){\line(1,-1){15}}

\put(-28,36){\line(0,-1){18}}\put(5,36){\line(0,-1){18}}

\put(-26,16){\line(1,0){28}}\put(33,57){\line(-1,-1){15}}\put(35,57){\line(1,-1){15}}

\put(17,36){\line(0,-1){18}}\put(50,36){\line(0,-1){18}}

\put(19,16){\line(1,0){28}}\put(36,61){\line(5,3){50}}

\put(19,40){\line(1,0){28}}\put(66,42){\line(-2,5){19}}\put(66,42){\line(2,5){19}}

\put(-12,90){\line(-1,-3){16}}\put(-11,90){\line(3,-2){43}}

\put(90,20){\line(1,0){39}}\put(108,38){\line(-5,-4){20}}\put(112,38){\line(5,-4){20}}

\put(110,42){\line(-5,2){120}}\put(88,22){\line(-3,5){40}}\put(131,22){\line(-2,3){46}}

\put(23,130){\line(1,0){24}}\put(20,133){\line(1,1){12}}\put(50,133){\line(-1,1){12}}

\put(20,127){\line(1,-1){12}}\put(50,127){\line(-1,-1){12}}

\put(18,128){\line(-3,-4){28}}\put(52,128){\line(4,-5){30}}

\put(-80,160){\line(1,0){20}}\put(-70,140){\vector(0,1){20}}\put(-70,120){\vector(0,-1){20}}\put(-80,100){\line(1,0){20}}

\put(-80,125){$C(G)$}

\put(-80,100){\line(1,0){20}}\put(-70,95){\vector(0,1){5}}\put(-70,85){\vector(0,-1){5}}\put(-80,80){\line(1,0){20}}

\put(-80,85){$A(G)$}

\put(-80,0){\line(1,0){20}}\put(-70,50){\vector(0,1){30}}\put(-70,30){\vector(0,-1){30}}

\put(-80,35){$D(G)$}

\put(-12,150){\circle{5}}\put(83,150){\circle{5}}
\put(-11,90){\circle*{5}}\put(84,90){\circle*{5}}\put(48,90){\circle*{5}}
\put(-11,60){\circle{5}}\put(-28,39){\circle{5}}\put(5,39){\circle{5}}
\put(-28,15){\circle{5}}\put(5,15){\circle{5}}

\put(34,60){\circle{5}}\put(17,39){\circle{5}}\put(50,39){\circle{5}}
\put(17,15){\circle{5}}\put(50,15){\circle{5}}

\put(66,39){\circle{5}}\put(88,19){\circle{5}}\put(132,19){\circle{5}}\put(110,40)
{\circle{5}}

\put(20,130){\circle{5}}\put(50,130){\circle{5}}\put(35,115){\circle{5}}\put(35,145){\circle{5}}

\put(-135,-30){Figure\ 1 \ \ The Gallai-Edmonds decomposition of a
graph $G$}
\end{picture}
\end{center}
\end{figure}

\vspace{30mm}

We say that a graph $G$ is {\it star-uniform}  if all star-factors
of $G$ have the same number of components. This concept is motivated
by an open problem posed by Hartnell and Rall \cite{har}, which
asked to characterize the family of graphs that all its star-factors
have the same number of edges and is motivated by the minimum cost
spanning tree and the optimal assignment problems. Hartnell and Rall
characterized star-uniform graphs with girth at least five. Wu and
Yu \cite{wu} settled the case of star-uniform graphs with girth
three but without leaves.

    To prove that a graph $G$ is not star-uniform, we usually show
that $G$ contains two star-factors with different numbers of
components. In this paper, we use the factor-criticality and
domination number to characterize all star-uniform graphs with the
minimum degree at least two. Our proof is heavily relied on
Gallai-Edmonds Structure Theorem.

\section{Main Results}

Let $\mathcal{S}$ be a star-factor with the maximum number of
components among all star-factors of $G$. If we choose one edge from
each component of $\mathcal{S}$, it yields a matching $M$.
Conversely, suppose $M$ is a maximum matching in $G$, then $G-V(M)$
is an independent set, and for each edge $uv$ in $M$, $u$ and $v$
can not be adjacent to distinct vertices of $G-V(M)$ due to the
maximality  of $M$. For each isolated vertex $x$ in $G-V(M)$, we add
an edge $e\in E(G)$ joining $x$ to a vertex in $V(M)$ and obtain a
star-factor with $|M|$ components. Hence we have the following
proposition.

\begin{prop}\label{pro1}
Let $G$ be a connected graph. Then the maximum number of components
of star-factors in $G$ is equal to the number of edges of a maximum
matching in $G$ (e.g., the matching number).
\end{prop}

Proposition \ref{pro1} shows the relationship between the maximum
number of components of star-factors and the matching number.
Similarly, we have a proposition to relate the minimum number of
components of star-factors and the domination number.

\begin{prop}\label{pro2}
Let $G$ be a connected graph. Then the minimum number of components
of star-factors in $G$ is equal to the domination number
$\gamma(G)$.
\end{prop}

\noindent {\bf Proof.} Let $\mathcal{S}$ be a star-factor in $G$
with the minimum number of components. Then all centers in
$\mathcal{S}$ form a dominating set, so $\gamma(G)$ is not greater
than the minimum number of components of star-factors in $G$.

    Conversely, suppose $D$ is a dominating set of the minimum order.
Then every vertex of $V(G)-D$ has at least one neighbor in $D$ and
every vertex of $D$ has at least one neighbor in $V(G)-D$ since $G$
has no isolated vertices. Now we construct a bipartite graph $B$
with bipartition $(V(G)-D)\cup D$ and edge set $E(B)=\{uv\ | \ u\in
V(G)-D, v\in D \mbox{ and} \ uv \in E(G)\}$. Then $B$ has a
star-factor, which can be regarded as a star-factor of $G$. Since
the number of components of a star factor of $B$ is at most $|D|$,
it follows that $\gamma(G)=|D|$ is greater than or equal to the
minimum number of components of star-factors in $G$. Therefore the
proposition is proved.  \hfill$\Box$

\medskip
The following result of Ore \cite{ore} provided a bound for
domination numbers of graphs without isolated vertices is a
corollary of Proposition~\ref{pro2}.

\begin{thm}\label{or}
{\rm{(Ore \cite{ore})}} If a graph $G$ has no isolated vertex, then
$\gamma(G) \leq \lfloor \frac{V(G)}{2} \rfloor$.
\end{thm}

  Every star factor $\mathcal{S}$ of a bipartite graph $G$
with bipartition $ X\cup Y$ such that $|X|\leq |Y|$ has at most
$|X|$ components. So, if $\gamma(G)=|X|$, then, by
Proposition~\ref{pro2}, each component of $\mathcal{S}$ contains
exactly one vertex of $X$, and thus we can obtain a matching from
$\mathcal{S}$ that covers $X$. Therefore the following result
follows.

\begin{prop}\label{pro3}
Let $G$ be a connected bipartite graph with bipartition $X\cup Y$
such that $|X|\leq |Y|$. If $\gamma(G)=|X|$, then $G$ contains a
matching of $|X|$ edges.
\end{prop}

Combining Propositions \ref{pro1} and \ref{pro2}, the following
theorem is true for all star-uniform graphs.

\begin{thm}\label{thm3}
A connected graph $G$ is star-uniform if and only if the size of a
maximum matching of $G$ is equal to the domination number
$\gamma(G)$.
\end{thm}

Next is the main result of this paper.

\begin{thm}\label{thm2}
A connected graph $G$ with $\delta(G) \geq 2$ is star-uniform if and
only if $G$ is one of the graphs shown in Figure $2$ or a bipartite
graph with bipartition $X\cup Y$ such that $g(G)=4$ and
$\gamma(G)=|X|\le |Y|$.
\end{thm}

\begin{figure}[h,t]
\setlength{\unitlength}{1pt}
\begin{center}
\begin{picture}(40,100)(0,40)
\put(-195,150){\line(0,-1){35}}\put(-195,150){\line(1,0){70}}
\put(-195,115){\line(1,0){70}}\put(-125,150){\line(1,-1){18}}
\put(-125,115){\line(1,1){18}}

\put(-195,150){\circle* {3}}\put(-160,150){\circle* {3}}
\put(-195,115){\circle* {3}}\put(-160,115){\circle* {3}}
\put(-125,150){\circle* {3}}\put(-125,115){\circle* {3}}
\put(-107,132){\circle* {3}}

\put(-85,150){\line(0,-1){35}}\put(-85,150){\line(1,0){70}}
\put(-85,115){\line(1,0){70}}\put(-15,150){\line(1,-1){18}}
\put(-15,115){\line(1,1){18}}\put(-50,150){\line(0,-1){35}}

\put(-85,150){\circle* {3}}\put(-50,150){\circle* {3}}
\put(-85,115){\circle* {3}}\put(-50,115){\circle* {3}}
\put(-15,150){\circle* {3}}\put(-15,115){\circle* {3}}
\put(3,132){\circle* {3}}

\put(25,150){\line(0,-1){35}}\put(25,150){\line(1,0){70}}
\put(25,115){\line(1,0){70}}\put(95,150){\line(1,-1){18}}
\put(95,115){\line(1,1){18}}
\put(60,150){\line(1,-1){35}}\put(60,115){\line(1,1){35}}

\put(25,150){\circle* {3}}\put(60,150){\circle* {3}}
\put(25,115){\circle* {3}}\put(60,115){\circle* {3}}
\put(95,150){\circle* {3}}\put(95,115){\circle* {3}}
\put(113,132){\circle* {3}}

\put(135,150){\line(0,-1){35}}\put(135,150){\line(1,0){70}}
\put(135,115){\line(1,0){70}}\put(205,150){\line(1,-1){18}}
\put(205,115){\line(1,1){18}}\put(170,115){\line(1,1){35}}
\put(170,150){\line(1,-1){35}}\put(170,150){\line(0,-1){35}}

\put(135,150){\circle* {3}}\put(170,150){\circle* {3}}
\put(135,115){\circle* {3}}\put(170,115){\circle* {3}}
\put(205,150){\circle* {3}}\put(205,115){\circle* {3}}
\put(223,132){\circle* {3}}

\put(-195,80){\line(0,-1){35}}\put(-195,80){\line(1,0){70}}
\put(-195,45){\line(1,0){70}}\put(-125,80){\line(1,-1){18}}
\put(-125,45){\line(1,1){18}}\put(-160,80){\line(1,-1){35}}
\put(-160,80){\line(0,-1){35}}

\put(-195,80){\circle* {3}}\put(-160,80){\circle* {3}}
\put(-195,45){\circle* {3}}\put(-160,45){\circle* {3}}
\put(-125,80){\circle* {3}}\put(-125,45){\circle* {3}}
\put(-107,62){\circle* {3}}

\put(-80,80){\line(0,-1){35}}\put(-80,80){\line(1,0){35}}
\put(-80,45){\line(1,0){35}}\put(-45,80){\line(1,-1){18}}
\put(-45,45){\line(1,1){18}}

\put(-80,80){\circle* {3}}\put(-80,45){\circle* {3}}
\put(-45,80){\circle* {3}}\put(-45,45){\circle* {3}}
\put(-27,62){\circle* {3}}

\put(10,80){\line(0,-1){35}}\put(10,80){\line(1,0){35}}
\put(10,45){\line(1,0){35}}\put(45,80){\line(1,-1){18}}
\put(45,45){\line(1,1){18}}\put(45,80){\line(0,-1){35}}

\put(10,80){\circle* {3}}\put(10,45){\circle* {3}}
\put(45,80){\circle* {3}}\put(45,45){\circle* {3}}
\put(63,62){\circle* {3}}

\put(100,80){\line(0,-1){35}}\put(100,80){\line(1,0){35}}
\put(100,45){\line(1,0){35}}\put(135,80){\line(1,-1){18}}
\put(135,45){\line(1,1){18}}\put(135,80){\line(0,-1){35}}
\put(100,45){\line(3,1){52}}

\put(100,80){\circle* {3}}\put(100,45){\circle* {3}}
\put(135,80){\circle* {3}}\put(135,45){\circle* {3}}
\put(153,62){\circle* {3}}

\put(210,70){\line(-1,-1){20}}\put(210,70){\line(1,-1){20}}
\put(190,50){\line(1,0){40}}

\put(210,70){\circle* {3}}\put(190,50){\circle* {3}}
\put(230,50){\circle* {3}}

\put(-5,15){Figure\ 2}
\end{picture}
\end{center}
\end{figure}

\vspace{10mm}

  Hence, there are only nine star-uniform graphs containing odd cycles
but there are infinite bipartite star-uniform graphs. Note that the
infinite family of bipartite graphs mentioned in Theorem~\ref{thm2}
satisfies $2=|X|\leq|Y|$ or $3\leq |X|<|Y|$.\\

\noindent {\bf Proof of Theorem \ref{thm2}.} It is easy to check
that all graphs shown in Figure~1 are star-uniform.  Let $G$ be a
bipartite graph with bipartition $X\cup Y$ such that
$\gamma(G)=|X|\le |Y|$. Then $G$ is star-uniform by
Proposition~\ref{pro3} and Theorem~\ref{thm3}.
\medskip

    Conversely, we shall prove that any connected star-uniform
graph with the minimum degree at least two is one of graphs given in
the theorem. Let $G$ be a connected star-uniform graph with
$\delta(G) \geq 2$. We consider the following cases.

\medskip
{\it Case 1.~ $G$ has a perfect matching.}
\medskip

  Let $M$ be a perfect matching of $G$. Since $G$ has a
perfect matching, $|V(G)|$ is even. If $|V(G)|=4$, then $G$ can only
be a $4$-cycle as $G$ is star-uniform. If $|V(G)|>4$, then there
exist three distinct edges $\{v_{1}v_{2},v_{3}v_{4},v_{5}v_{6}\}$ of
$M$ such that $v_2v_3, v_4v_5 \in E(G)$ since $G$ is connected and
$\delta(G) \geq 2$. Then we can find  two stars $T_{1}$ and $T_{2}$
with centers $v_{2}$ and $v_{5}$, respectively, which cover
$\{v_{1}v_{2},v_{3}v_{4},v_{5}v_{6}\}$. Thus
$\{T_{1},T_{2},M-\{v_{1}v_{2},v_{3}v_{4},v_{5}v_{6}\}\}$ is a
star-factor of $G$ with $|M|-1$ components, but $M$ is a star-factor
with $|M|$ components, a contradiction to star-uniform of $G$.

    Therefore $G$ is a 4-cycle, which is one of the bipartite graphs
given in the theorem.

\medskip
{\it Case 2.~ $G$ has no perfect matching.}
\medskip

  Let $M$ be a maximum matching of $G$. By Gallai-Edmonds
Structure Theorem, we know that each component of the subgraph
induced by $D(G)$ is factor-critical and the subgraph induced by
$C(G)$ has a perfect matching. Moreover, $M$ consists of a
near-perfect matching of each component of $D(G)$, a perfect
matching of each component of $C(G)$ and a matching which matches
all vertices of $A(G)$ to vertices in distinct components of $D(G)$.
Since $G$ is connected,  for each component $D$ of $D(G)$, there is
at least one vertex in $D$ which is adjacent to a vertex in $A(G)$.
As $D$ is factor-critical, without loss of generality, we may assume
that each isolated vertex in $V(G)-V(M)$ is adjacent to a vertex in
$A(G)$. Now we add an edge $e\in E(G)$ between each isolated vertex
in $V(G)-V(M)$ and some vertex in $A(G)$, then we obtain a special
star-factor $\mathcal{S}$ with $|M|$ components such that all
vertices in $A(G)$ are centers of stars in $\mathcal{S}$.

    In the following discussion, we often delete some edges
from $\mathcal{S}$ and then add other edges in $E(G)-\mathcal{S}$ to
$\mathcal{S}$ to construct another star-factor with different number
of components and thus yields a contradiction to that $G$ is
star-uniform.

\medskip
{\it Claim 1.~  $C(G)=\emptyset$.}
\medskip

Let $C$ be a component of $C(G)$. Since $G$ is connected, then there
exists an edge $u_{1}v_{1}\in \mathcal{S}$ in $C$ such that $u_{1}$
is adjacent to a vertex $x$ in $A(G)$, and $v_{1}$ is also adjacent
to a vertex in $C$ or $A(G)$ since $\delta(G)\geq 2$. If $v_{1}$ is
adjacent to a vertex $y$ in $A(G)$, then deleting edge $u_{1}v_{1}$
from $\mathcal{S}$ and adding two edges $u_{1}x$ and $v_{1}y$ to
$\mathcal{S}$, we obtain another star-factor  with $|M|-1$
components. If $v_{1}$ is adjacent to a vertex $v_{2}$ in $C$, then
deleting edge $u_{1}v_{1}$ from $\mathcal{S}$ and adding two edges
$u_{1}x$ and $v_{1}v_{2}$ to $\mathcal{S}$, we also obtain a
star-factor with $|M|-1$ components. Either case yields a
contradiction.

\medskip
{\it Claim 2.~ $A(G)$ is an independent set.}
\medskip

Suppose $uv$ is an edge in the subgraph induced by $A(G)$ and
$T_{u}$ is a star in $\mathcal{S}$ with center $u$. For each leaf
$x$ in $T_{u}$, where $x$ is in a component $D$ of $D(G)$, we
perform the following operation: if $D$ is singleton, then $x$ is
adjacent to another vertex $y$ in $A(G)$ (since $\delta(G)\geq 2$)
and we delete the edge $ux$ from $\mathcal{S}$ and add the edge $xy$
to $\mathcal{S}$; if $D$ is nontrivial, then $x$ is adjacent to a
vertex $z$ in $D$, we delete the edge $ux$ from $\mathcal{S}$ and
add the edge $xz$ to $\mathcal{S}$. By adding the edge $uv$ to
$\mathcal{S}$, we obtain another star-factor $\mathcal{S}'$ with
$|M|-1$ components since $v$ is a center of some star in
$\mathcal{S}$, a contradiction.

\medskip
{\it Claim 3.~ $A(G)=\emptyset$ or all components of $D(G)$ are
singletons.}
\medskip

Suppose $A(G)\neq \emptyset$ and a component $D$ in $D(G)$ is
nontrivial. Then there exists a star $T_{v}$ in $\mathcal{S}$ such
that a leaf $x$ of $T_v$ is contained in $D$ and the center $v$ is
in $A(G)$. If $T_v$ has another leaf $x'$, we delete the edge $vx'$
from $\mathcal{S}$ and add another edge $x't$ to $\mathcal{S}$,
where $t\in A(G)\cup D(G)$ is the center of a star of $\mathcal{S}$.
Now we obtain another star-factor with the same number of components
as $\mathcal{S}$. So, without loss of generality, we may assume that
$v$ has no other leaves except $x$ in $T_v$. Suppose $e_{1},e_{2},
\cdots, e_{m}$ are $K_{1,1}$-stars of $\mathcal{S}$ in $D$. Then $x$
is adjacent to a vertex in $D$, but can not be adjacent to the two
vertices of a certain star $e_{i}$ $(1\leq i\leq m)$, otherwise the
four vertices in $xv$ and $e_{i}$ can be covered by one star with
center $x$, a contradiction. So let $e_{i}=yz$ be a star with only
one vertex, say $y$,  adjacent to $x$ and $z$ is adjacent to a
vertex $u \in D\cup A(G)$ since $\delta(G)\geq 2$. Then by removing
$yz$ from $\mathcal{S}$, and adding $yx$ and $zu$ to $\mathcal{S}$,
we can obtain another star-factor with $|M|-1$ components, a
contradiction.

\medskip
{\it Subcase 2.1.} $A(G)=\emptyset$.
\medskip

Since $G$ is connected, $A(G)=\emptyset$ implies that $G$ is
factor-critical.

\medskip
{\it Claim 4.} If $G$ is factor-critical, then $|V(G)|\leq 7$.
\medskip

  Since $G$ is factor-critical, then for each edge $xy$, both $G-x$
and $G-y$ have perfect matchings, denoted by $M_{x}$ and $M_{y}$,
respectively. Let $H=M_{x} \cup M_{y}$.  Then $H$ is a spanning
subgraph of $G$ and contains an alternating path $P_{xy}$ connecting
$x$ and $y$. Without loss of generality, we may assume that all
components except $P_{xy}$ are independent multiple edges and the
alternating path $P_{xy}$ is a longest path among all pairs
$(M_x,M_y)$ of perfect matchings in $G-x$ and $G-y$ over all edges
in $E(G)$. The order of $P_{xy}$, denoted by $p$, is odd. In the
following, we assume $|V(G)|\geq 9$ and then construct two
star-factors of $G$ with different numbers of components, and thus
yields a contradiction.

If $p\geq 9$, then $P_{xy}\cup \{xy\}$ has two star-factors
${\mathcal{S}_{1}}$ and ${\mathcal{S}_{2}}$ with
$\lfloor\frac{p}{2}\rfloor$ stars and $\lfloor\frac{p}{2}\rfloor-1$
stars, respectively. So $H\cup \{xy\}$ contains two star-factors
$\mathcal{S}_{1}\cup(M_{x}-P_{xy})$ and
$\mathcal{S}_{2}\cup(M_{x}-P_{xy})$ with different numbers of
components, a contradiction.

If $p=7$, then there is at least one vertex $u$ in $P_{xy}$ which is
adjacent to a vertex $v_{1}$ in $H-P_{xy}$ since $|V(G)|\geq 9$ and
$G$ is connected. Then $P_{xy}\cup \{xy,uv_{1},v_{1}v_{2}\}$ has two
star-factors $\mathcal{S}_{3}$ and $\mathcal{S}_{4}$ with three
stars and four stars, respectively. That is, $G$ contains two
star-factors $\mathcal{S}_{3}\cup(M_{x}-P_{xy}-v_{1}v_{2})$ and
$\mathcal{S}_{4}\cup(M_{x}-P_{xy}-v_{1}v_{2})$ with different
numbers of components.

If $p=5$, then there exists an edge $v_{1}v_{2}$ in $H-P_{xy}$ which
is joined to $P_{xy}$ by an edges $v_{1}u_{1}$, where $u_{1}\in
V(P_{xy})$ since $G$ is connected and $|V(G)|\geq 9$.  Suppose that
$v_{2}$ is adjacent to another edge $v_{3}v_{4}$ in
$H-P_{xy}-v_{1}v_{2}$. Suppose that $v_{2}$ is adjacent to $v_{3}$
in $G$, then $P_{xy}\cup
\{xy,v_{1}v_{2},v_{3}v_{4},u_{1}v_{1},v_{2}v_{3}\}$ has two
star-factors $\mathcal{S}_{5}$ and $\mathcal{S}_{6}$ with three
stars and four stars, respectively. So $\mathcal{S}_{5}\cup
\{M_{x}-P_{xy}-v_{1}v_{2}-v_{3}v_{4}\}$ and $\mathcal{S}_{6}\cup
\{M_{x}-P_{xy}-v_{1}v_{2}-v_{3}v_{4}\}$ are two star-factors with
different numbers of components in $G$.  Thus $v_{2}$ is not
adjacent to any vertex in $H-P_{xy}$. In this case, $v_{2}$ must be
adjacent to a vertex of $P_{xy}$ as $\delta(G) \ge 2$.  If $v_{1}$
and $v_{2}$ have a common neighbor $u_{2}$, then $P_{xy}\cup
\{xy,v_{1}v_{2},v_{1}u_{2},v_{2}u_{2}\}$ can be decomposed into two
stars or three stars, so $G$ contains two star-factors with
different numbers of components. If $v_{1}$ and $v_{2}$ are adjacent
to two distinct vertices $u_{3}$ and $u_{4}$ in $P_{xy}$
respectively, then $u_{3}$ and $u_{4}$ are nonadjacent in
$P_{xy}\cup \{xy\}$ since $P_{xy}$ is a longest alternating path.
Then $P_{xy}\cup \{xy,v_{1}v_{2},v_{1}u_{3},v_{2}u_{4}\}$ can be
decomposed into two stars or three stars, we again obtain two
star-factors with different numbers of components in $G$.

If $p=3$, then there exists an edge $v_{1}v_{2}$ in $H-P_{xy}$ such
that $v_1u$ is an edge of $G$, where $u\in P_{xy}$. If  $v_{2}$ is
incident with another edge $v_{3}v_{4}$ in $H-P_{xy}-v_{1}v_{2}$,
where $v_{2}v_3$ is an edge of $G$, then $P_{xy}\cup
\{xy,v_{1}v_{2},v_{3}v_{4},uv_{1},v_{2}v_{3}\}$ can be decomposed
into two stars or three stars. So $H\cup \{xy,uv_{1},v_{2}v_{3}\}$
has two star-factors with the different numbers of components.
Otherwise both $v_{1}$ and $v_{2}$ are only adjacent to vertices in
$P_{xy}$.  If $v_2$ is adjacent to $u'\neq u$ in $G$, denote
$z=P_{xy}-u-u'$, then $M_{v_{1}}=M_x-P_{xy}-v_{1}v_{2}+v_{2}u'+uz$
and $M_{v_{2}}=M_x-P_{xy}-v_{1}v_{2}+v_{1}u+u'z$ are maximum
matching of $G-v_1$ and $G-v_2$, respectively, and the path
connecting $v_1$ and $v_2$ in $M_{v_{1}}\cup M_{v_{2}}$ is of length
$5$, which is a contradiction to the choice of $P_{xy}$. Hence
$v_{2}$ is also adjacent to $u$ in $G$.  Then  $P_{xy}\cup
\{xy,v_{1}v_{2},uv_{1},uv_{2}\}$ can be decomposed into one star or
two stars, we obtain two star-factors of $G$ with different numbers
of components. Consequently Claim~4 is proved.

\medskip
All factor-critical graphs of order three, five or seven, with
$\delta(G)\geq 2$ and $\gamma(G)= \lfloor \frac{V(G)}{2} \rfloor$
are shown in Figure~2.

\medskip
 {\it Subcase 2.2.~~ $A(G)\ne \emptyset$ and all components of
$D(G)$ are singletons.}

  By Claim~1, this assumption implies that $G$ is a bipartite graph
with bipartition $A(G)\cup D(G)$.

\medskip
{\it Claim 5.~~ $g(G)=4$ and $\gamma(G)=|A(G)|$.}
\medskip

By Gallai-Edmonds Structure Theorem, $|A(G)|<|D(G)|$ and every
maximum matching of $G$ covers all vertices in $A(G)$. So
$\gamma(G)=|A(G)|$ by Theorem \ref{thm3}.

Since $G$ is a bipartite graph, $g(G)=4,6,8,\cdots$. Suppose that
$g(G)\ge 6$. Let $C$ be  a girth cycle of $G$. If a vertex not
contained in $C$ is adjacent to two distinct vertices of $C$, then
there is a cycle shorter than $g(C)$. Hence any vertex in
$V(G)-V(C)$ is adjacent to at most one vertex of $V(C)$. Thus by
deleting all the edges incident with the cycle $C$ but not on it, we
obtain a spanning subgraph $H$ without isolated vertices. But $C$ is
a component of $H$ and $C$ can be decomposed into  $g(G)/2$ stars or
$(g(G)/2)-1$ stars. This is a contradiction since $H$ is a spanning
subgraph of $G$. Hence $g(G)=4$. Since  $G$ has no perfect matching,
we have $2\leq |A(G)|<|D(G)|$.

 This completes the proof of Theorem~\ref{thm2}. \hfill$\Box$\\

\medskip

    Next we construct an infinite family of star-uniform
connected bipartite graphs. A connected graph that has no cut
vertices is called a {\it block}. A {\it block of a graph} is a
subgraph that is a block and is maximal with respect to this
property: Every graph is the union of its blocks. Let $G$ be a
connected graph such that each block of $G$ is $K_{2,2}$ and each
block has a pair of nonadjacent vertices of degree two in $G$, e.g.,
the graph shown in Figure $3$ is such a graph. Then $G$ is
star-uniform. To see this, let $u$ and $v$ be two nonadjacent
vertices of degree two in a block $B$, then neither $u$ nor $v$ can
be a center of a star $K_{1,n}$ $(n\geq 2)$ in any star-factor of
$G$. For a star $K_{1,1}$, any its vertex can be considered as a
center, so we designate other two vertices rather than $u$ and $v$
in $B$ as the centers in any star-factor of $G$. Hence all
star-factors of $G$ have the same number of components, i.e., $G$ is
star-uniform.

\begin{figure}[h,t]
\setlength{\unitlength}{1pt}
\begin{center}
\begin{picture}(40,100)(0,60)

\put(-60,150){\line(1,0){30}}\put(-60,150){\line(1,-1){20}}
\put(-40,130){\line(1,0){30}}\put(-30,150){\line(1,-1){20}}

\put(-60,150){\circle* {3}}\put(-40,130){\circle* {3}}
\put(-30,150){\circle* {3}}\put(-10,130){\circle* {3}}

\put(-10,130){\line(1,0){30}}\put(-10,130){\line(1,-1){20}}
\put(10,110){\line(1,0){30}}\put(20,130){\line(1,-1){20}}

\put(-10,130){\circle* {3}}\put(10,110){\circle* {3}}
\put(20,130){\circle* {3}}\put(40,110){\circle* {3}}

\put(40,110){\line(1,0){30}}\put(40,110){\line(1,-1){20}}
\put(60,90){\line(1,0){30}}\put(70,110){\line(1,-1){20}}

\put(40,110){\circle* {3}}\put(60,90){\circle* {3}}
\put(70,110){\circle* {3}}\put(90,90){\circle* {3}}

\put(40,110){\line(0,1){20}}\put(40,110){\line(2,1){20}}
\put(60,120){\line(0,1){20}}\put(40,130){\line(2,1){20}}

\put(60,120){\circle* {3}}\put(40,130){\circle* {3}}
\put(60,140){\circle* {3}}

\put(40,110){\line(0,-1){20}}\put(40,110){\line(-2,-1){20}}
\put(20,100){\line(0,-1){20}}\put(20,80){\line(2,1){20}}

\put(20,100){\circle* {3}}\put(40,90){\circle* {3}}
\put(20,80){\circle* {3}}

\put(20,80){\line(0,-1){20}}\put(20,80){\line(-2,-1){20}}
\put(0,70){\line(0,-1){20}}\put(0,50){\line(2,1){20}}

\put(0,70){\circle* {3}}\put(20,60){\circle* {3}}
\put(0,50){\circle* {3}}

\put(-10,130){\line(0,-1){20}}\put(-10,130){\line(-2,-1){20}}
\put(-30,120){\line(0,-1){20}}\put(-30,100){\line(2,1){20}}

\put(-30,120){\circle* {3}}\put(-10,110){\circle* {3}}
\put(-30,100){\circle* {3}}

\put(-5,20){Figure\ 3}
\end{picture}
\end{center}
\end{figure}

\vspace{10mm}

An {\it edge-weighting} of a graph $G$ is a function $w:
E(G)\longrightarrow \mathbb{N}^+$, where $\mathbb{N}^+$ is a set of
positive integers. The {\it weight} of a star-factor $\mathcal{S}$
in $G$ under $w$ is the sum of all the weight values for edges
belonging to $\mathcal{S}$, i.e., $w(\mathcal{S})=\Sigma_{e\in
E(\mathcal{S})}w(e)$. Now it is nature to ask the following question
which is proposed in \cite{har} and still open.

\vspace{2mm}

\begin{oq}
For a given graph $G$, does there exist an edge-weighting $w$ of $G$
such that every star-factor of $G$ has the same weights under $w$?
\end{oq}

%\vskip 20pt
%
%\noindent \title{\Large\bf Acknowledgments} \maketitle
% The authors are indebted to the anonymous referees for their constructive suggestions.


\begin{thebibliography}{99}

\bibitem{amahashi} A. Amahashi and M. Kano, On factors with
given components, {\em Discrete Math.},
42(1982), pp. 1-6.

\bibitem{bb} B. Bollob\'{a}s, {\it Modern Graph Theory}, 2nd Edition,
Springer-Verlag New York, Inc. 1998.

\bibitem{ed} J. Edmonds, Paths, trees, and flowers, {\em Canad. J. Math.},
17(1965), pp. 449-467.

\bibitem{har} B. L. Hartnell and D. F. Rall,
On graphs having uniform size star
factors, {\em Australas. J. Combin.},
34(2006), pp. 305-311.

\bibitem{hay} T. W. Haynes, S. T. Hedetniemi and P. J. Slater,
{\it Fundamentals of Domination in Graphs}, Marcel Dekker, Inc., New
York, 1998.

\bibitem{lm} L. Lov\'{a}sz and M.D. Plummer,
{\it Matching Theory}, North-Holland Inc., Amsterdam, 1986.

\bibitem{ore} O. Ore, {\it Theory of Graphs},
AMS Publication 38, Providence, RI, 1962.

\bibitem{wu} Y. Wu and Q. Yu, Uniform Star-factors of Graphs with Girth
Three, (submitted).

\bibitem{yu} Q. Yu, Counting the number of star-factors in graphs,
{\em J. Combin. Math. Combin. Comput.},
23(1997), pp. 65-76.

\end{thebibliography}
\end{document}